\documentclass[reqno, 12pt]{amsart}

\def\ssign{\textsection\nobreak\hspace{1pt plus 0.3pt}}
\makeatletter
\let\origsection=\section 
\def\mysection{\@mystartsection{section}{1}\z@{.7\linespacing\@plus\linespacing}{.5\linespacing}{\normalfont\scshape\centering\ssign}}
\def\section{\@ifstar{\origsection*}{\mysection}}
\def\appendix{\par\c@section\z@ \c@subsection\z@
   \let\sectionname\appendixname
   \let\section=\origsection
   \def\thesection{\@Alph\c@section}}
\def\@mystartsection#1#2#3#4#5#6{\if@noskipsec \leavevmode \fi
 \par \@tempskipa #4\relax
 \@afterindenttrue
 \ifdim \@tempskipa <\z@ \@tempskipa -\@tempskipa \@afterindentfalse\fi
 \if@nobreak \everypar{}\else
     \addpenalty\@secpenalty\addvspace\@tempskipa\fi
 \@dblarg{\@mysect{#1}{#2}{#3}{#4}{#5}{#6}}}
\def\@mysect#1#2#3#4#5#6[#7]#8{\edef\@toclevel{\ifnum#2=\@m 0\else\number#2\fi}\ifnum #2>\c@secnumdepth \let\@secnumber\@empty
  \else \@xp\let\@xp\@secnumber\csname the#1\endcsname\fi
  \@tempskipa #5\relax
  \ifnum #2>\c@secnumdepth
    \let\@svsec\@empty
  \else
    \refstepcounter{#1}\edef\@secnumpunct{\ifdim\@tempskipa>\z@ \@ifnotempty{#8}{\@nx\enspace}\else
        \@ifempty{#8}{.}{\@nx\enspace}\fi
    }\@ifempty{#8}{\ifnum #2=\tw@ \def\@secnumfont{\bfseries}\fi}{}\protected@edef\@svsec{\ifnum#2<\@m
        \@ifundefined{#1name}{}{\ignorespaces\csname #1name\endcsname\space
        }\fi
      \@seccntformat{#1}}\fi
  \ifdim \@tempskipa>\z@ \begingroup #6\relax
    \@hangfrom{\hskip #3\relax\@svsec}{\interlinepenalty\@M #8\par}\endgroup
    \ifnum#2>\@m \else \@tocwrite{#1}{#8}\fi
  \else
  \def\@svsechd{#6\hskip #3\@svsec
    \@ifnotempty{#8}{\ignorespaces#8\unskip
       \@addpunct.}\ifnum#2>\@m \else \@tocwrite{#1}{#8}\fi
  }\fi
  \global\@nobreaktrue
  \@xsect{#5}}
\makeatother

\usepackage[T1]{fontenc}
\usepackage{lmodern}
\usepackage[babel]{microtype}
\usepackage[english]{babel}

\usepackage{geometry}
\usepackage{mathdots}
\usepackage[
sorting=nty,
backend=biber,
natbib=true,
bibstyle=ieee,
citestyle=numeric-comp,
isbn=false,
maxbibnames=9,
backref=true,
uniquename=init
]{biblatex}

\DefineBibliographyStrings{english}{  backrefpage={$\upharpoonleft$}, 
backrefpages={$\upharpoonleft$}}
\usepackage{csquotes}
\DeclareFieldFormat{title}{\mkbibquote{#1}}
\DeclareFieldFormat[misc]{date}{Preprint (#1)}

\renewbibmacro*{doi+eprint+url}{
	\iftoggle{bbx:url}
	{\iffieldundef{doi}{\usebibmacro{url+urldate}}{}}{}
	\newunit\newblock
	\iftoggle{bbx:eprint}
	{\usebibmacro{eprint}}
	{}
	\newunit\newblock
	\iftoggle{bbx:doi}
	{\printfield{doi}}
	{}
}

\renewbibmacro{in:}{}

\DeclareFieldFormat{page}{#1}

\usepackage{csquotes}

\usepackage{mathrsfs}
\usepackage[colorlinks=true, linkcolor=blue, filecolor=magenta, urlcolor=cyan]{hyperref}
\usepackage{xcolor}
\hypersetup{
    colorlinks,
    linkcolor={red!80!black},
    citecolor={red!60!black},
    urlcolor={blue!60!black}
 }

\usepackage{amsmath,amssymb,amsthm}
\usepackage{mathrsfs}
\usepackage{mathabx}\changenotsign
\usepackage{dsfont}

\numberwithin{equation}{section}
\numberwithin{figure}{section}

\usepackage{enumitem}
\newcommand{\customlabel}[2]{#2\def\@currentlabel{#2}\label{#1}}
\def\rmlabel{\upshape({\itshape \roman*\,})}

\def\alabel{\upshape({\itshape \alph*\,})}

\usepackage{tikz}
\usepackage{scalerel}

\newsavebox\vdegbox
\savebox\vdegbox{\tikz{
		\draw[black,fill=black] (90:1) circle (.35);
		\draw[black,line width=0.10cm] (210:1) circle (.30);
		\draw[black,line width=0.10cm] (330:1) circle (.30);
		\draw[opacity=0] (0:1.2) circle (0.1);
	}}

\newsavebox\vvbox
\savebox\vvbox{\tikz{
		\draw[black,line width=0.10cm] (90:1) circle (.30);
		\draw[black,fill=black] (210:1) circle (.35);
		\draw[black,fill=black] (330:1) circle (.35);
		\draw[opacity=0] (0:1.2) circle (0.1);
	}}

\newsavebox\pdegbox
\savebox\pdegbox{\tikz{
		\draw[black,line width=0.10cm] (90:1) circle (.30);
		\draw[black,fill=black] (210:1) circle (.35);
		\draw[black,fill=black] (330:1) circle (.35);
		\draw[black,line width=0.28cm ] (210:1) -- (330:1);
		\draw[opacity=0] (0:1.2) circle (0.1);
	}}

\newsavebox\vvvbox
\savebox\vvvbox{\tikz{
		\draw[black,fill=black] (90:1) circle (.35);
		\draw[black,fill=black] (210:1) circle (.35);
		\draw[black,fill=black] (330:1) circle (.35);
		\draw[opacity=0] (0:1.2) circle (0.1);
	}}

\newsavebox\evbox
\savebox\evbox{\tikz{
		\draw[black,fill=black] (90:1) circle (.35);
		\draw[black,fill=black] (210:1) circle (.35);
		\draw[black,fill=black] (330:1) circle (.35);
		\draw[black,line width=0.28cm ] (210:1) -- (330:1);
		\draw[opacity=0] (0:1.2) circle (0.1);
	}}

\newsavebox\eebox
\savebox\eebox{\tikz{
		\draw[black,fill=black] (90:1) circle (.35);
		\draw[black,fill=black] (210:1) circle (.35);
		\draw[black,fill=black] (330:1) circle (.35);
		\draw[black,line width=0.28cm ] (90:1) -- (330:1);
		\draw[black,line width=0.28cm ] (90:1) -- (210:1);
		\draw[opacity=0] (0:1.2) circle (0.1);
	}}

\newsavebox\eeebox
\savebox\eeebox{\tikz{
		\draw[black,fill=black] (90:1) circle (.35);
		\draw[black,fill=black] (210:1) circle (.35);
		\draw[black,fill=black] (330:1) circle (.35);
		\draw[black,line width=0.28cm ] (90:1) -- (330:1);
		\draw[black,line width=0.28cm ] (90:1) -- (210:1);
		\draw[black,line width=0.28cm ] (210:1) -- (330:1);
		\draw[opacity=0] (0:1.2) circle (0.1);
	}}

\theoremstyle{plain}
\newtheorem{thm}{Theorem}[section]

\newtheorem{clm}{Claim}

\theoremstyle{definition}

\newtheorem{exmp}[thm]{Example}

\newenvironment{claimproof}[1][{Proof}]{\begin{proof}[#1]}{\end{proof}}

\let\eps=\varepsilon
\let\theta=\vartheta
\let\rho=\varrho
\let\phi=\varphi

\DeclareMathOperator{\ex}{ex}

\let\polishlcross=\l
\def\l{\ifmmode\ell\else\polishlcross\fi}

\makeatletter
\def\moverlay{\mathpalette\mov@rlay}
\def\mov@rlay#1#2{\leavevmode\vtop{%
   \baselineskip\z@skip \lineskiplimit-\maxdimen
   \ialign{\hfil$\m@th#1##$\hfil\cr#2\crcr}}}
\newcommand{\charfusion}[3][\mathord]{
    #1{\ifx#1\mathop\vphantom{#2}\fi
        \mathpalette\mov@rlay{#2\cr#3}
      }
    \ifx#1\mathop\expandafter\displaylimits\fi}
\makeatother

\def\qqand{\qquad\text{and}\qquad}

\newcommand{\vrhup}[1]{\scaleobj{0.6}{\scalerel*{\rightharpoonup}{#1}}}
\newcommand{\nrhup}{\mathord{\scaleobj{0.6}{\scalerel*{\rightharpoonup}{x}}}}

\def\vseq#1{\ThisStyle{%
  \mathord{\vbox{\offinterlineskip\ialign{%
    \hfil##\hfil\cr
    $\SavedStyle{}_{\smash{\vrhup#1}}$\cr
    \noalign{\kern-0.7\scriptspace}
    $\SavedStyle#1$\cr}}}}}

\def\seq#1{\ThisStyle{%
  \mathord{\vbox{\offinterlineskip\ialign{%
    \hfil##\hfil\cr
    $\SavedStyle{}_{\smash{\nrhup}}$\cr
    \noalign{\kern-0.5\scriptspace}
    $\SavedStyle#1$\cr}}}}}

\let\setminus=\smallsetminus

\let\to=\lra

\usepackage{xspace}

\makeatletter
\newcommand{\pushright}[1]{\ifmeasuring@#1\else\omit\hfill$\displaystyle#1$\fi\ignorespaces}
\newcommand{\pushleft}[1]{\ifmeasuring@#1\else\omit$\displaystyle#1$\hfill\fi\ignorespaces}
\makeatother

\title{The codegree Tur\'an density of $3$-uniform tight cycles}

\author[S. Piga]{Sim\'on Piga}
\address[S. Piga]{Fachbereich Mathematik, Universit\"at Hamburg, Hamburg, Germany}
\email{simon.piga@uni-hamburg.de }
\author[N. Sanhueza-Matamala]{Nicol\'as Sanhueza-Matamala}
\address[N. Sanhueza-Matamala]{Departamento de Ingeniería Matemática, Facultad de Ciencias Físicas
y Matemáticas, Universidad de Concepción, Chile}
\email{nicolas@sanhueza.net}
\author[M.~Schacht]{Mathias Schacht}
\address[M.~Schacht]{Fachbereich Mathematik, Universit\"at Hamburg, Hamburg, Germany}
\email{schacht@math.uni-hamburg.de}

\thanks{}


\begin{document}

\begin{abstract}
    Given any~$\eps>0$ we prove that every sufficiently large~$n$-vertex~$3$-graph~$H$ where every pair of vertices is contained in at least~$(1/3+\eps)n$ edges contains a copy of~$C_{10}$, i.e.\ the tight cycle on~$10$ vertices. 
    In fact we obtain the same conclusion for every cycle~$C_\ell$ with~$\ell\geq 19$.
\end{abstract} 

\maketitle

\section{Introduction}
We consider an extremal problem for hypergraphs. A $k$-uniform hypergraph $H$ is defined by a vertex set~$V(H)$ and a set of edges~$E(H)\subseteq V(H)^{(k)} = \{S \subseteq V(H)\colon \vert S \vert = k\}$. 
Throughout this note, unless specified otherwise, we refer to $3$-uniform hypergraphs simply as hypergraphs. 
For a given hypergraph~$F$, the \textit{extremal number~$\ex(n,F)$} for~$n$ vertices is the maximum number of edges in an~$n$-vertex hypergraph that does not contain a copy of~$F$. 
The \textit{Tur\'an density} $\pi(F)$ is defined as
$$ \lim_{n \to \infty} \frac{\ex(n,F)}{\binom{n}{3}}\,,$$ 
this is well-defined for every~$F$, since the sequence~${\ex(n,F)}/{\binom{n}{3}}\geq 0$ is non-increasing.
Determining the Tur\'an densities of hypergraphs is a central problem in combinatorics. 
Despite considerable efforts by many researchers, Tur\'an densities are known only for few hypergraphs.
For discussion of techniques, results, and variations, see the surveys by Keevash~\cite{Keevash}, Balogh, Clemen, and Lidick\'y~\cite{balogh} and Reiher~\cite{Christian}.

Our focus here is on the variation called \textit{codegree Tur\'an density}, introduced by Mubayi and Zhao~\cite{MZ}. 
Given a hypergraph~$H$ and a subset~$S \in V^{(2)}$, the \textit{neighbourhood}~$N_H(S)$ and \textit{codegree}~$d_H(S)$ of~$S$ are defined by
$$N_H(S) = \{v\in V(H) \colon S\cup v\in E(H)\} 
\qqand 
d_H(S) = |N_H(S)|\,,$$
and when~$H$ is clear from the context we will omit it from the notation.
We will omit unnecessary parenthesis and commas from the set-theoretic notation, and in particular we write~$d(uv)$ instead of~$d(\{u,v\})$. 
The minimum codegree of~$H$ among all possible sets~$S$ of size two is denoted by~$\delta_2(H)$. 
For a hypergraph~$F$ and an integer~$n$, the \emph{codegree Tur\'an number}~$\ex_2(n,F)$ is the maximum~$d$ such that there exists an~$F$-free hypergraph~$H$ on~$n$ vertices with~$\delta_2(H) \geq d$; and the \emph{codegree Tur\'an density of $F$} is
\[
\gamma(F) = \lim_{n \to \infty} \frac{\ex_2(n,F)}{n}\,,
\]
which is also well-defined for every $F$ \cite[Proposition 1.2]{MZ}.
Clearly,~$\gamma(F) \leq \pi(F)$.

Similarly, the codegree Tur\'an density is known for only few hypergraphs (see, e.g.,~\cite[Table~1]{balogh}).
In particular, a computer-assisted proof by Falgas-Ravry, Pikhurko, Vaughan, and Volec~\cite{K4-flag} determined that~$\gamma(K_4^{-}) = 1/4$, where $K_4^-$ is the hypergraph obtained from~$K_4$ by removing one edge. 
In contrast, $\gamma(K_4)$ is not known, with Czygrinow and Nagle~\cite{codegK4} conjecturing that~$\gamma(K_4) = 1/2$.

Given two hypergraphs $F, G$, an \emph{homomorphism} from $F$ to $G$ is a map $\phi\colon V(F) \to V(G)$ which preserves edges, i.e. $\phi(e) \in E(G)$ for each $e \in E(H)$.
By the phenomenon of supersaturation (see, e.g.~\cite[Section 2]{Keevash} and \cite[Proposition 1.4]{MZ}) it turns out that if there exists an homomorphism from $F$ to $G$ then $\pi(F) \leq \pi(G)$ and $\gamma(F) \leq \gamma(G)$.

We are interested in studying codegree Turán densities of tight cycles.
Given an integer~$\ell \geq 3$, a \emph{tight cycle~$C_\ell$} is a hypergraph with vertex set~$\{v_1, \ldots, v_\ell\}$ and edge set~$\{v_iv_{i+1}v_{i+2}\colon i \in \mathbb{Z}/\ell\mathds{Z}\}$. 
Whenever $\ell$ is divisible by $3$ we have that $C_\ell$ is $3$-partite, which implies that $\gamma(C_\ell) = \pi(C_\ell) = 0$ (see \cite{erdos64}), so the interesting cases concern $\ell$ not divisible by~$3$ only.
Recently, Kam\v{c}ev, Letzter, and Pokrovsky~\cite{KamcevLetzterPokrovsky2024} determined $\pi(C_\ell)$ for those values of~$\ell$, as long as $\ell$ is sufficiently large. 

The following lower bound construction shows that $\gamma(C_\ell) \geq 1/3$ for~$\ell$ not divisible by $3$.

\begin{exmp}\label{example}
    Let~$n\in 3\mathds N$ and let~$H=(V,E)$ be a hypergraph where $V=V_1\dot\cup V_2\dot\cup V_3$ with~$|V_i|=n/3$ and 
    $$E=\{uvw \in V^{(3)}\colon u\in V_i, v\in V_j, w\in V_k \text{ and } i+j+k \equiv 1\bmod 3\}\,.$$
It is easy to check that~$\delta_2(H)\geq n/3-1$. 
Let $C_\ell$ be a cycle in $H$, we will show that~$\ell$ is divisible by $3$.
For each $v \in V(C_\ell)$ let $c(v) = j$ if $v \in V_j$ and
set $\Omega = \sum_{e \in E(C_\ell)} \sum_{v \in e} c(v)$.
By construction, we have that $\sum_{v \in e} c(v) \equiv 1 \bmod 3$ for each $e \in E(C_\ell)$, therefore~$\Omega \equiv \ell \bmod 3$.
Moreover, since every vertex of $C_\ell$ is contained exactly in $3$ edges, we also have that $\Omega \equiv 0 \bmod 3$.
Hence, $\ell \equiv 0 \bmod 3$.
\end{exmp}

The previously known best upper bound for codegree Turán densities of tight cycles is due to Balogh, Clemen, and Lidick\'y~\cite{balogh}.
One of their results yields~$\gamma(C_\ell)\leq 0.3993$ for every~$\ell\geq 5$ except~$\ell=7$.

In this note we establish an upper bound matching Example~\ref{example} for almost every~$\ell$ not divisible by three. 

\begin{thm} \label{thm:main}
    For $\ell \in \{10, 13, 16\}$ and for every $\ell \geq 19$ not divisible by $3$, $\gamma(C_\ell) = 1/3$.
\end{thm}



We can use homomorphisms to obtain codegree Tur\'an densities for longer cycles using shorter ones.
Indeed, note that for each $\ell \geq 6$, there is an homomorphism from $C_{\ell + 3}$ to $C_\ell$ (by wrapping around the last three vertices), and therefore $\gamma(C_{\ell + 3}) \leq \gamma(C_{\ell})$.
Moreover, for any $t \geq 2$, there is an homomorphism from $C_{t \ell}$ to $C_\ell$ (by transversing $C_\ell$ $t$ times), so $\gamma(C_{t\ell}) \leq \gamma(C_{\ell})$ also holds.
Combining these two observations, it is easy to see that we only need to prove~$\gamma(C_{10})\leq 1/3$.


\section{Proof of Theorem~\ref{thm:main}}

Given~$\eps>0$ let~$n_0\in \mathds N$ be sufficiently large and let~$H$ be a hypergraph on~$n\geq n_0$ vertices with~$\delta_2(H)\geq (1/3+\eps)n$. 
It suffices to show that $H$ contains an homomorphic image of $C_{10}$.
For a contradiction, suppose not. 
We separate the rest of the proof in a series of claims. 

\begin{clm}\label{clm:manyK4-}
    Every edge of $H$ is contained in a copy of~$K_4^{-}$.
\end{clm}
\begin{claimproof}
    Let~$e=xyz\in E(H)$ and note that~$d(xy)+d(xz)+d(yz) \geq (1+3\eps)n$. 
    Hence, there is a vertex~$v\in V\setminus \{x,y,z\}$ such that~$v$ is in two neighbourhoods~$N(xy), N(xz), N(yz)$.
    Suppose~$v\in N(xy)\cap N(xz)$.
    Then the edges $\{ xyz, xyv, xzv\}$ form a copy of $K_4^-$.
\end{claimproof}

We say the only vertex of degree $3$ in a $K_4^{-}$ is the \emph{apex} of $K_4^{-}$.
We say that a pair of distinct vertices~$u,v\in V(H)$, is an \textit{apex pair} if there is a copy of~$K_4^{-}$  containing~$u$ and $v$, where either~$u$ or $v$ is the apex.
Similarly, we say they are a~\textit{base pair} if there is a copy of~$K_4^{-}$ containing~$u$ and~$v$ where neither of them is the apex.

\begin{clm}\label{clm:redblue}
    Every pair of distinct vertices is either an apex pair or a base pair, but not both. 
\end{clm}
\begin{claimproof}
    Observe that Claim~\ref{clm:manyK4-} together with the minimum codegree condition imply that every pair of vertices is contained in a copy of~$K_4^{-}$.
    In particular, every pair is an apex pair or a base pair.

    Suppose that the pair~$uv$ is simultaneously an apex pair and a base pair.
    Consequently, we can assume that there are~$K$ and~$K'$, copies of~$K_4^{-}$, both containing the vertices~$u$ and~$v$ and such that~$v$ is the apex of~$K$ and neither~$u$ nor $v$ is the apex of~$K'$. 
    Let~$V(K)=\{u,v, x,y\}$ and~$V(K')=\{u,v,a,b\}$ be the vertex sets of~$K$ and~$K'$ respectively, where $a$ is the apex of~$K'$. 
    Observe that the ordering~$(\boldsymbol v,u,\boldsymbol{a},b,v,\boldsymbol a,u,\boldsymbol v,x,y)$ forms an homomorphic copy of~$C_{10}$, where we marked the apexes for clarity. 
\end{claimproof}


We define an auxiliary directed graph~$D$ with on the vertex set~$V(H)$ with arcs given by
$$E(D)=\big\{(u,v)\in V(H)^2 \colon \text{$uv$ is an apex pair with apex $v$}\big\}\,.$$

\begin{clm}\label{clm:oriented}
    $D$ does not contain a directed cycle of length $2$.
\end{clm}
\begin{claimproof}
    Suppose ~$(a,x),(x,a)\in E(D)$.
    Then there are~$K$ and~$K'$, copies of~$K_4^{-}$, both containing the vertices~$a$ and~$x$, and such that~$a$ is the apex of~$K$ and $x$ is apex of~$K'$. 
    Let~$V(K)=\{a,x, b,c\}$ and~$V(K')=\{a,x,y,z\}$ be the vertex sets of~$K$ and~$K'$ respectively. 
    Observe that the ordering~$(\boldsymbol x,\boldsymbol a,y,\boldsymbol x,z,\boldsymbol a,\boldsymbol x,b,\boldsymbol a,c)$ forms an isomorphic copy of~$C_{10}$, where we marked the apexes for clarity. 
\end{claimproof}

Let~$B = \{uv\in V(H)^{(2)}\colon \text{ $uv$ is a base pair}\}$ and note that due to Claims~\ref{clm:redblue} and~\ref{clm:oriented} for every pair~$uv\in V(H)^{(2)}$ exactly one of the following alternatives hold:
\begin{enumerate}[label=\rmlabel]
    \item \label{it1} $(u,v)\in E(D)$, 
    \item \label{it2}$(v,u)\in E(D)$, or
    \item \label{it3} $uv\in B$.
\end{enumerate}

The following claim shows that the edges of $B$ and the arcs of~$D$ are strongly related. 
\begin{clm}\label{clm:degrees}
    For every $v\in V(H)$ we have:
    \begin{enumerate}[label=\alabel]
        \item \label{it:Bimplies+} If~$d_B(v) >0$, then $d^+_D(v)\geq (1/3+\eps)n$.
        \item If~$d^{+}_D(v)>0$, then~$d_B(v)\geq (1/3+\eps)n$.
        \item If~$d^-_D(v)>0$, then~$d^-_D(v)\geq (1/3+\eps)n$.
    \end{enumerate}
\end{clm}
\begin{claimproof}
    Since the proofs are all analogous, we only show~\ref{it:Bimplies+}. 
    Let~$u$ be such that~$uv\in B$ and let~$w\in N(uv)$ chosen arbitrarily. 
    Due to Claim~\ref{clm:manyK4-} there is a~$K_4^{-}$ containing the edge~$uvw$.
    Observe that neither $u$ nor $v$ can be the apex of such~$K_4^{-}$, otherwise, $uv$ would be an apex pair, contradicting Claim~\ref{clm:redblue}. 
    Therefore~$w$ is the apex, and~$(v,w)\in E(D)$. 
    Hence~$N(uv)\subseteq N_D^+(v)$, meaning that~$|N_D^+(v)|\geq (1/3+\eps)n$.
\end{claimproof}

Finally, if there is an vertex~$v^\star\in V(H)$ with~$d_B(v^\star) >0$, $d^{+}_D(v^\star)>0$, and $d^-_D(v^\star)>0$, then Claim~\ref{clm:degrees} yields a contradiction with Claim~\ref{clm:redblue} or Claim~\ref{clm:oriented}, since there would be a pair for which two of the three alternatives~\ref{it1}, \ref{it2}, \ref{it3} hold. 
We shall find such vertex~$v^\star$. 

First, suppose there are two distinct vertices~$u,v$ with~$d_D^{+}(u)=d_D^+(v)=0$. 
Then~$uv\in B$, due to Claim~\ref{clm:redblue}, and in particular~$d_B(u),d_B(v)>0$. 
However, Claim~\ref{clm:degrees} yields a contradiction, since this implies~$d_D^+(u),d_D^+(v)>0$. 
Hence, there is at most one vertex having zero out-degree in~$D$.

Secondly, take two disjoint edges~$e_1$ and~$e_2$ and note that Claim~\ref{clm:manyK4-} implies that there are vertices~$v_1\in e_1$ and~$v_2\in e_2$ with~$d_D^{-}(v_1),d_D^{-}(v_2)>0$. 
One of them, say~$v_1$, has positive out-degree as well, i.e.\ $d^{+}_D(v_1)>0$. 
Since Claim~\ref{clm:degrees} yields~$d_B(v_1)>0$ we are done by taking~$v^{\star} = v_1$.

\section{Concluding remarks}

It would be interesting to settle the remaining values of~$\gamma(C_\ell)$. 
The case $\ell = 4$ is equivalent to the determination of $\gamma(K_4)$ and, as mentioned in the introduction, Czygrinow and Nagle~\cite{codegK4} conjectured~$\gamma(K_4) = 1/2$.
It seems plausible that Example~\ref{example} is optimal for all other values of~$\ell$ not divisible by three.
In other words, that~$\gamma(C_\ell)=1/3$ for every~$\ell\geq 5$ not divisible by three. 
Note that by our previous remarks, for this result it would suffice to show~$\gamma(C_5)\leq 1/3$ and  $\gamma(C_7) \leq 1/3$.

Determining whether Example~\ref{example} is optimal for~$\ex(n,C_\ell)$ for~$\ell\geq 5$ not divisible by three is a natural question. 
We believe a more careful analysis of the proof of Theorem~\ref{thm:main} yields a constant~$c\in \mathds N$ such that $$\ex_2(n,C_{10}) \leq \frac{n}{3} + c\,,$$
for sufficiently large~$n$ and finding the optimal constant~$c$ remains open.

Finally, for~$k$-uniform hypergraphs with~$k\geq 4$, the problem of determining~$\gamma\big(C_\ell^{(k)}\big)$ in general remains open. 
For general lower-bound constructions see~\cite[Section 10]{HanLoSanhueza}.

\printbibliography

\end{document}